\documentclass[10pt]{article}
\usepackage{amsmath,amssymb,amsthm,graphics,graphicx}
\usepackage[T1]{fontenc}
\usepackage{multicol}
\usepackage{epstopdf}
\begin{document}
\vspace{5pt}
\begin{center}
\vspace{2pt}
    {\Large{Prime Power Divisibility, Periodicity, and Other Properties of Some Second Order  Recurrences}
 }
\\
{Soumyabrata Pal \footnote{Fourth Year student at IIT Kharagpur}   \\}
{Shankar M. Venkatesan \footnote {Bangalore, India}}
 \end{center}
 \vspace{4pt}
\begin{abstract}
Wall published a paper in 1960 on the Fibonacci sequence where he derived many results concerning the period and prime power divisibility modulo m. His periodicity results have been generalized before to second order linear recurrences.  Here we study the sequences generated by such recurrences, with starting values of \{0,1\}: Among other things, we (a) derive new prime power divisibility results not only for the general recurrence, but also interestingly for the squares of such sequences, (b) derive the period by new methods, (c) establish many novel identities modulo a prime power, (d) show derivations involving powers of matrices generated by these general recurrences, etc.
\end {abstract}

\section{Introduction}
Fibonacci sequences have been extensively studied in the past. While this prior work is highly interesting, we shall describe in this introduction only what is relevant to the main themes of this paper, which includes the extension of the prime power divisibility and other results of Wall. After Wall [3], Robinson [4] derived  periodicity using the Fibonacci matrix method (for a history of this method see [18]). General recurrences and over finite fields were studied by many including [5, 6]. Then came a breakthrough result of Matijasevich [8] that uses a Fibonacci identity (extended to general recurrrences in Theorem 9). On a diffferent note, Knuth and Wilf [11] studied the interesting prime power divisibility of generalized binomial coefficients generalizing the well-known work of Kummer. The power of the logarithmic bounding implied by Wall's prime power result was used in [12] to show that the Fibonacci sequences written in any base is not context-free (our result in the current paper generalizes this result as well in a strong sense). General recurrences were also presented in a nice exposition in[19] and shown to share many of the well known basic properties of Fibonacci and Lucas sequences. Periodicity was studied in [16] in the context of a splitting field characterization. Prime power divisibility modulo m was further studied in [13,14,15, 17] and many related results were proved, but none of these extend the prime power results of Wall to general recurrences.

Here we study sequences generated by general second order linear recurrences (with starting values of 0 and 1), defined by the equation e(n)=Ae(n-1)+Be(n-2) where A and B are integers. For our known Fibonacci sequence A and B are 1 with starting values 0 and 1. Among other things, we generalize the prime power divisibility results of Wall on Fibonacci recurrences modulo m, show new derivations involving powers of matrices generated by these recurrences, prove many novel identities, study the starting values, etc.

\section{Arithmetic Index Progressions in Recurrences}
Wall proved in 1960 [3, Theorem 3] that the indices of elements in the Fibonacci sequence that are zero modulo m form an arithmetic progression. We will extend Wall's Method to obtain the same result for general recurrences (this will be needed later).\\
Theorem 1: For the general sequence $e(n)=Ae(n-1)+Be(n-2)$ the terms $ e(n)=0 (mod ~m)$ forms a simple arithmetic progression except in the case for which $gcd(B,m) \neq 1$.\\ 
Proof:  We prove by induction that the following identity (well known for the Fibonacci sequence since B=1 for Fibonacci sequence) holds for the general recurrence 
$e(n+t)=e(n+1)e(t)+Be(n)e(t-1)$ (this like many other results is equally easy to prove with fibonacci matrices [19], but we reserve this for attacking some significant problems later).\\
 Obviously for $e(0)=0$ and $e(1)=1$ and $e(2)=b$ the relation trivially holds true for $t=1$ and $t=2$. Say it holds for all $t <=k$. Now for $t=k+1$ we have that \\
$e(n+1)e(k+1)+Be(n)e(k)$\\
$=e(n+1)[Ae(k)+Be(k-1)]+e(n)[B*Ae(k-1)+B*Be(k-2)]$\\
$=A[e(n+1)e(k)+B*e(n)e(k-1)]+B[e(n+1)e(k-1)+B*e(n)e(k-2)]$\\
$=Ae(n+k)+Be(n+k-1)$\\
$=e(n+k+1)$\\
Hence  the relation is true for e(n) using Wall's [3] arguments on setting $e(i)=0(mod ~ m)$ and $e(j)=0(mod  ~m)$ and finding that $e(i+j)=0(mod ~ m)$. But for the case where $gcd(B,m) \neq 1$, then the sequence goes into a cycle where e(n) does not equal 0 for all n.\\
This is because for the period to exist we need e(n)=0 and e(n+1)=1 for some n. Now it means that B.e(n-1)=1(mod  ~m) or $e(n-1)=B^{-1}$(mod  ~m) but $B^{-1}  (mod ~m)$ does not exist.Hence the period does not exist. We can in fact state that 0 will also not appear, otherwise one can go backwards if gcd(a,m)=1 and we get back 0 and 1 but this is also not possible since it goes into a cycle where 0 and 1 do not exist consecutively. But the fact that 0 will exist on the periodic cycle of $e(n)(mod  ~m)$ depends on the fact that one can go backwards in the cycle . So in this case the lemma does not hold.

\section{Prime Power Divisibility in General Sequences}
Here u(n) is taken to be our normal Fibonacci sequence and v(n) is taken to be the Lucas sequence with starting values 2 and 1. Wall [3] studied prime power repetition and divisibility in Fibonacci Sequences :\\

Wall's Theorem: If $u(n)=0 (mod  ~p^{e})$,  but $ \neq 0(mod ~ p^{e+1})$ is the first occurrence for r such that $u(r)=0(mod ~ p^{e})$ then $u(pn)=0 (mod ~ p^{e+1} )$  but $\neq 0 (mod ~ p^{e+2})$ is the first occcurence for $u(x)=0 (mod  ~p^{e+1})$. \\
Proof:According to Binet's formula $u(n )=(r^{n}-s^{n})/ \sqrt(5)$ where r and s are conjugate surds obtained by solving for the n'th Fibonacci number as a function of n and $v(n)=r^{n}+s^{n}$. Also $u(an)=(r^{an}-s^{an})/$$\sqrt(5)$ . So solving for $r^{an}$ and $s^{an}$ and using the binomial theorem and the property that v(n)=v(n-1)+v(n+1), we get that\\
$gcd(v(n),e(n))=1$ or $2$ and\\
$u(an)=2^{1-a}*(^{a}C_{j})_{j=odd}5^{(j-1)/2}u(n)^{j}v(n)^{a-j}$\\
So $u(an)=2^{1-a}u(n)(Ku(n)^{2}+av(n)^{a-1})$ .\\
Hence for u(an) to be divisible by $p^{e+1}$ p must divide a given that $p^{e}$ divides u(n).\\

We now extend Wall's theorem above to the general recurrences.\\
Theorem 2: In sequence e (E) defined before, if $e(n)=0 (mod  ~p^{k})$,  but $ \neq 0(mod ~ p^{k+1})$ is the first occurrence for r such that $e(r)=0(mod ~ p^{k})$ then $e(pn)=0 (mod ~ p^{k+1} )$  but $\neq 0 (mod  ~p^{k+2})$ is the first occcurence for $x$ such that $e(x)=0 (mod ~ p^{k+1})$. \\
Proof: For $e(n)=Ae(n-1)+Be(n-2)$, by Binet's formula, we have
$e(n)=k(r^{n}-s^{n})$ where $r+s=A$ and $rs=-B$ (unless specified otherwise we will assume $e(0)=0$ and $e(1)=1$).
We need to prove the existence of a sequence $v(n)$ such that $v(n)=k'(r^{n}+s^{n})$ for a constant $k'$.
Taking up from here, we can define $v(n)=pe(n-1)+qe(n+1)$. We know that $v(n)=Av(n-1)+Bv(n-2)$. Also we have $e(n+1)=Ae(n)+Be(n-1)$.
Solving these 2 equations we have
Bv(n)=e(n+1)(pB+q)-A.q.e(n).\\
Assuming that $gcd(A,B)=1$ , we reach the conclusion that $B$ must divide $q$ for $v(n)$ to be an integer for all n. So let $q=xB$ where $x$ is an integer. Substituting, we have\\
$v(n)=(p+x)e(n+1)-Axe(n)$\\
$v(n)=k[(p+x)(r^{n+1}-s^{n+1})-Ax(r^{n}-s^{n})]$\\
$v(n)=k[r^{n}((p+x)r-Ax)+s^{n}(Ax-(p+x)s))]$\\
Now the coefficients of $r^{n}$ and $s^{n}$ must be equal. So we get \\
$(p+x)(r+s)=2Ax$ or $(p+x)A=2Ax$ or $p=x$\\
or $p=q/B$

Hence the sequence exists $v(n)$ exists if we choose $p$ and $q$ such that $q=pB$. So choosing $p=1$, $q=B$ we have\\ $v(n)=2e(n+1)-Ae(n)$.\\
Now since $gcd(A,B)=1$ we must have that $gcd(e(n),e(n+1))=1$ . Hence we must have that $gcd(v(n),e(n))=1 ~or~ 2$.\\ 
We have $e(an)=k(r^{an}-s^{an}) $\\
$r^{n}=(1/2)*(e(n)/k+v(n)/k')$ and $s^{n}=(1/2)*(-e(n)/k+v(n)/k')$\\
So, $e(an)=(k/(2)^{a})[(e(n)/k+v(n)/k')^{a}-(v(n)/k'-e(n)/k)^{a}]$\\ 
Now, $k'=k*[(p+x)r-Ax]$ .\\
 Since $p=x=1$,$ k'=k*(2r-A)$ or $ k'=k(r-s)$ or $ k'=1$. Since By Binet's formula $k=1/(r-s)$  $k'$must be $1$.\\
Hence,  $e(an)=(k/(2)^{a})[(e(n)/k+v(n))^{a}-(v(n)-e(n)/k)^{a}]$ Also 1/k must be of the form $\sqrt b$ because r and s are the roots of an equation whose coefficients are real.\\
$e(an)=2^{1-a}*(^{a}C_{j})_{j=odd}k^{(1-j)}e(n)^{j}v(n)^{a-j}$ .\\
 So $e(an)=2^{1-a}e(n)(Ke(n)^{2}+av(n)^{a-1})$ .\\
 This formula along with the fact that  $gcd(v(n),e(n))=1 ~or~ 2$ along with the result of section on arithmetic progressions leads to the conclusion that the law of repetion of primes also holds for the general recurrence relation with the same argument as in the case of Fibonacci sequence. However we must remember that this proof is only for starting values of 0 and 1.  \\

While it is obvious then that there is always a member E(n) in the general sequence which has more trailing zeros (when written in base m) than any given k, we also get the following from the above theorem (please see [12]):\\ 
Corollary: The number of trailing zeros of any member E(n) in the general sequence written in any base m,  is $<= k. log (n)$ for some fixed k.\\
This exponential gap can be a powerful tool in studying the recognition complexity of these and other sequences. While Automata Theory is not an object of study in this paper, we observe that this gap (even a non-linear one will suffice) is sufficient to prove that every such generalized sequence written in base m is non-context-free using the same technique as in [12]. 

\section{Alternative Analysis using Fibonacci Matrices}
This section tries to analyze an approach using the Fibonacci matrices culminating in a very interesting theorem when this approach is matched with that of the previous one which follows Wall's work.
Let us define the Fibonacci matrix\\
$E_k$=
$\begin{pmatrix}
e(k+1) & e(k)\\
e(k) & e(k-1)\\
\end{pmatrix}$
and the matrix A as
$\begin{pmatrix}
1 & 1\\
1 & 0\\
\end{pmatrix}$
\\
So we have that 
$E_n=A^{n}*E_0$
where $E_0= 
\begin{pmatrix}
e(1) & e(0)\\
e(0) & e(-1)\\
\end{pmatrix}$
\\
So we have again 
$E_{na}=A^{na}*E_{0}$\\
or $E_{na}=(E_{n}*(E_{1})^{-1})^{a}*E_{1}$\\
or $E_{na}=(E_{n})^{a}*(E_{1})^{1-a}$\\
Now taking determinant on both sides we observe that
$\begin{vmatrix}
E_{0}\end{vmatrix}
=1$\\
given that e(0)=0 and e(1)=1 and e(-1)=1 since e() is the fibonacci sequence in this case.
Hence, we have that\\
$e(na+1)e(na-1)-e(na)^2=(e(n+1)e(n-1)-e(n)^2)^a*(1)^{1-a}$\\
or $ e(na+1)e(na-1)-e(na)^2=(e(n+1)e(n-1))^{a}-e(n)^{2a}-a*e(n)^{2}*K$ where K contains the rest of the terms.
Hence, we have that\\
$e(na+1)e(na-1)-e(na)^2=(e(n+1)e(n-1)-e(n)^2)^a*(1)^{1-a}$\\
or $ e(na+1)e(na-1)-e(na)^2=(e(n+1)e(n-1))^{a}-e(n)^{2a}-a*e(n)^{2}*K$\\
So assuming that we have e(n) is the first number for which $e(n)=0$ (mod $p$) but $\neq 1 (mod$ $p^{2})$. Also from the previous section we know that e(np) is the first term for which $e(.)$ is divisible by $p^{2}$\\
On taking both sides congruent to $p^{2}$ and putting a=p,\\
$e(np+1)e(np-1)-e(np)^2=e(n+1)e(n-1)^{p}-e(n)^{2p}+p*e(n)^{2}*K$,  for some K.
From the previous section we know that if $e(n)=0(mod$ $p$) then $p^{2}$ divides e(np). Using that fact we find that since $p$ divides $e(n)$ , $p^{3}$ divides $p*e(n)^{2}*K-e(n)^{2p}$. Since the $gcd(e(n),e(n+1))$ and $gcd(e(n),e(n-1)=1$ , we have the following novel result which has to be satisfied if the theorem in the previous section has to hold true.\\
Theorem 4: $e(np+1)e(np-1)=(e(n+1)e(n-1))^{p}(mod~$ $p^{2})$ if $e(n)=0(mod~$ $p$)\\
In that case, $p^{3}$ divides the left hand side as well. But since the left hand side is the square of an integer hence $p^{4}$ must divide the left hand side.Hence $p^{2}$ divides e(np) and e(np) is the first term which is divisible by $p^{2}$ but not $p^{3}$. So  in this case the 2 approaches coincide.\\ 

This idea can also be similarly extended to $p^{e}$ from $p$ by using the same arguments and hence the extension is\\
Theorem 5:  $e(np+1)e(np-1)=(e(n+1)e(n-1))^{p}(mod$ $p^{e+1})$ if $e(n)=0(mod$ $p^{e}$)\\

\section{Recurrence Spaces}
We can define the concept of a recurrence space induced by $\{A,B\}$, as also had been previously done in [19].\\
Lemma: The sequence $ e(n)=A.e(n-1)+B.e(n-2)$  forms a space for which any $1$ member forms the basis.\\
Proof:- Say $e(n)=A.e(n-1)+B.e(n-2)$ and we define $v(n)=R.e(p)+S.e(q)$ . \\
Then the sequence v obtained as below also belongs to the space. 
$ v(n)=R[Ae(p-1)+Be(p-2)]+S[Ae(q-1)+Be(q-2)]$\\
$=v(n)=A[Re(p-1)+Se(q-1)]+B[Re(p-2)+Se(q-2)]$\\
$=>v(n)=Av(n-1)+Bv(n-2)$\\
Corollary: This also holds for multiterm general recurrence relation i.e $v(n)=a_1e(p)+a_2e(q)+a_3e(r)+....$  then v(n) belongs to the class V which follows the property $e(n)=a_1e(n-1)+a_2e(n-2)+a_3e(n-3)+....$
An interesting point to note here will be that the recurrence space defined above can also be extended to Fibonacci Matrices. We have defined v(n) to be
v(n)=Re(p)+Se(q). \\
So writing it out in the matrix form using the defintion of Fibonacci matrices given in the previous section, we have\\
$V_n=RA^{p}E_0+SA^{q}E_0$\\
or $V_0=RA^{p-n}E_0+SA^{q-n}E_0$\\
or $V_0=BE_0$ where $B=RA^{p-n}+SA^{q-n}$\\
and we have $V_n=A^{n}BE_0$\\
or $V_n=A^{n}V_0$\\
So this implies that recurrence space is also applicable to the Fibonacci matrices.\\
This recurrence space concept will allow us to relate general sequences with different starting numbers to sequences with starting values as 0 and 1, as shown in few of the subsequent theorems \\

\section{Extension to the General Recurrence}
For the general recurrence , the technique used above to analyse the fibonacci matrices does not apply. The problem lies in the fact that e(-1)  depends on the number m since B.e(-1)=1(mod  m). So to extend this result to the general recurrence we have to take a slightly different path.\\
Again we can define,\\
$E_n=A^{n-1}*E_1$
where $E_1=
\begin{pmatrix}
e(2) & e(1)\\
e(1) & e(0)\\
\end{pmatrix}$
\\
and $A=
\begin{pmatrix}
a & b\\
1 & 0\\
\end{pmatrix}$
\\
Then we have that $E_{an}=A^{an-1}*E_1$
or $E_{an}*A=A^{an}*E_1$
or $E_{an}*A=(E_n*A*(E_1)^(-1))^{a}*E_1$

So taking determinant on both sides, we have \\
$\begin{vmatrix}
E_{1}\end{vmatrix}
=-1$ and
$\begin{vmatrix}
A\end{vmatrix}
=-b$ \\
thus yielding the following result:\\
Lemma:  $e(na+1)e(na-1)-e(np)^2=(-b)^{a-1}*(e(n+1)e(n-1)-e(n)^2)^a*(-1)^{1-a}$\\

Since we know from the previous sections that the theorem from Section 4 holds for the general recurrence as well, it gives rise to a very interesting result below because only if it is satisfied the above equation will be 0(mod $p^{e+1}$). This inference can be drawn from the previous 2  sections and the above equation.\\
Theorem 6:
 $(-b)^{p-1}*(e(n+1)e(n-1))^p=e(np+1)e(np-1)(mod$ $p^{e+1}$) where $e(n)=0(mod$ $p^{e}$).\\
\section{Periodicity : Extension of Wall's work to general recurrences}
Let us first review Wall's work on periodicity once more:\\
Wall's theorem on periodicity:  If $k(p^{2}) \neq k(p)$ (the well known Wall's question) then $k(p^{e})=p^{e-1}k(p)$ and if t is the largest number such that $k(p^{t})=k(p)$ then $k(p^{e})=p^{e-t}k(p)$ where e>t for the Fibonacci sequence\\
Proof:-In earlier section we extended Wall’s theorem for the Arithmetc progression of the zeroes in the fibonacci sequence modulo a particular prime power and it’s higher power to the general recurrence.\\
Taking it up from there,Wall then proved that\\
$u(an+1)=2^{-a}(K*u(n)^{2}+a*u(n)*v(n)^{a-1}+v(n)^{a})$\\
We now know that u(nx) are the terms which are equal to 0(mod $p^{e}$) and u(pnx)=0 are the terms which are equal to 0(mod $p^{e+1}$). So setting u(pnx+1)=1(mod $p^{e+1}$) to obtain $k(p^{e+1})=pnx$  the equation above gives \\
 $(v(nx)/2)=1(mod$ $p^{e+1})$.\\
So $v(nx)=2( mod$ $p^{e})$ and thus u(nx+1)=1 (mod $p^{e}$). So we have that $k(p^{e+1})=pk(p^{e})$ and Wall's question whether $k(p^{2})=k(p)$ remains unanswered\\

So we work along the same lines as above and prove the following :\\
Theorem 7: If $k(p^{2}) \neq k(p)$ then $k(p^{e})=p^{e-1}k(p)$ and if t is the largest number such that $k(p^{t})=k(p)$ then $k(p^{e})=p^{e-t}k(p)$ where e>t for the general sequence\\
Proof: From earlier sections, we have obtained $v(n)=r^{n}+s^{n}$  and we also have $u(n)=k(r^{n}-s^{n})$ where k=1/(r-s) and r and s are the 2 conjugate surds which helps to relate u(n) and v(n) as a function of n. THe exact values of r and s can be calculated by using generating functions or using matrix theory to find the eigen values of the transforming matrix. So again we have that,\\
$e(an+1)=k(r^{an+1}-s^{an+1})$\\
 $e(an+1)=(k/(2)^{a})[r*(e(n)/k+v(n))^{a}-s*(v(n)-e(n)/k)^{a}]$\\
Since it is exactly the same form as the one obtained by Wall, the same conclusions also follow i.e we have 
$e(an+1)=2^{-a}(K*u(n)^{2}+a*u(n)*v(n)^{a-1}+v(n)^{a})$ after similar simplification. 
Again setting e(pnx+1)=1(mod $p^{e+1})$)
we have that $(v(nx)/2)=2(mod$ $p^{e}$) adnd thus we have e(nx+1)=1(mod $p^{e}$) from v(n)=2e(n+1)-Ae(n).
So the same theorem applies for the general recurrence as well. However the question of whether $k(p^{2})=k(p)$ exists for some p still remains unanswered.

\section{Periodicity using Fibonacci matrices}

Periodicity is easy to prove modulo a prime $p$ since only $p^2$ adjacent pairs can exist in the sequence, but it is possible to show a much smaller period by using Fermat's little Theorem as has been done by many authors. Periods for Fibonacci matrices were proved by Robinson [4], and for general recurrences by [16] among others. In this section we are going to prove an interesting theorem on periodicity for the general recurrence using the Fibonacci matrices.\\
Theorem 8: For a sequence e(n)=ae(n-1)+be(n-2) the square numbers $e(0)^2$,$e(1)^{2}$,... generated by squaring the members of the general recurrence  satisfies the property that  $k(p^{e})=p^{e-1}*k(p^{e})$ if $k(p^{2}) \neq k(p)$  and if t is the highest number such that $k(p^{t})=k(p)$ then $k(p^{e})=p^{e-t}*k(p)$ where e>t\\
Proof: By definition of the Fibonacci matrix for the general recurrence, we have 
$E_{n+1}n=A^{n}*E_1$ and $E_{pn+1}n=A^{pn}*E_1$ 
where $E_1=
\begin{pmatrix}
e(2) & e(1)\\
e(1) & e(0)\\
\end{pmatrix}$
\\ 
and $A=
\begin{pmatrix}
A & B\\
1 & 0\\
\end{pmatrix}$
\\ 
So by substituting like before, we have that\\
$E_{pn+1}=(E_n)^{p}*(E_1)^{1-p}$\\
Taking determinant on both sides and observing that determinant of $E_1$=-1 and 1-p is an even number, we have that\\
$e(pn+2)e(pn)-e(pn+1)^{2}=(e(n+2)e(n)-e(n+1)^{2})^{p}$ \\
So if e(n) is the first term such that e(n)=0(mod $p^{e}$) then pn is the first term for which e(pn)=0(mod $p^{e+1}$). So if we take congruent to $p^{e+1}$ on both sides and set $ e(n+1)=1(mod$ $p^{e}$) so that $k(p^{e})=n$ we have ,by using binomial expansion that
$e(pn+1)^{2}=e(n+1)^{2p}(mod$ $p^{e+1})$ and hence $e(pn+1)^{2}=1(mod$ $p^{e+1}$)
This concludes the proof of the theorem since it means that  $k(p^{e+1})=pn=p^{e-1}*k(p^{e})$. Again the possibilty of $k(p^{2})=k(p)$ remains open.                                                                                                                                                                                                                                                                                                                                                                                                                                                                                                                                                                                                                                                                                                                                                                                                                                                                                                                                                                                                                                                                                                                                                                                                                                                                                                                                                                                                                                                                                                                                                                                                                                                                                                                                                                                                                                                                                                                                                                                                                                                                                                                                                                                                                                                                                                                                                                                                                                                                                                                                                                                                                                                                                                                                                                                                                                                                                                                                                                                                                                                                                                                                                                                                                                                                                                                                                                                                                                                                                                                                                                                                                                                                                                                                                                                                                                                                                                                                                                                                                                                                                                                                                                                                                                                                                                                                                                                                                                                                                                                                                                                                                                                                                                                                                                                                                                                                                                                                                                                                                                                                                                                                                                                                                                                                                                                                                                                                                                                                                                                                                                                                                                                                                                                                                                                                                                                                                                                                                                                                                                                                                                                                                                                                                                                                                                                                                                                                                                                                                                                                                                                                                                                                                                                                                                                                                                                                                                                                                                                                                                                                                        

\section{ Starting Points of Period in Generalized Recurrence} 
Going back to section1, we stated that 0 may not the starting point of the periodic cycle in e(n)=ae(n-1)+be(n-2) if $gcd(b,p) \neq 1$ .
So the interesting objective of this section is to investigate the nature of the starting point.
Say that for the series e(n)=ae(n-1)+be(n-2) (mod p), the first 3 terms are:
0,1,a. So, for 1,a to be in the periodic cycle we must have ,\\
a=a+bx (mod p)\\
or bx=0(mod p)\\
or x=tp/gcd(p,b) (mod p) for some integer $t\neq0$ since 0 cannot be in the cycle\\
So the next question is the existence of e(n-1) if  e(n)=x, e(n+1)=1 and e(n+2)=a.
For that we must have ,\\
$tp/gcd(p,b) +be(n-1)=1(mod$ $p$)\\
or be(n-1)=1-(tp/gcd(p,b))(mod p)\\
For the solution to exist, we must have that\\
gcd(b,p) | (1-tp/gcd(p,b))\\
Therefore,\\
tp/gcd(p,b)=1(mod gcd(p,b))\\
or $t=(p/gcd(p,b))^{-1}(mod$ gcd(p,b))\\
This is condition for t.
Hence e(0) should be\\
$[(p/gcd(p,b))^{-1}(mod$ $gcd(p,b)$)*($p/gcd(p,b))](mod$ $p$). \\
We can also further establish the starting point of a periodic cycle for any arbitrary v(0) and v(1) by making use of our results for recurrence spaces which states that v(n) lies in the recurrence space of E for which e(n) is the basis. So we can define v(n)=Re(n-1)+Se(n) and calculate R and S from the information given about v(n) and we can then derive the starting point of v(n) as well since we know the starting points of e(n).

\section{Related  Divisibility Results}

In section 4 we proved that  $e(an)=2^{1-a}e(n)(Ke(n)^{2}+av(n)^{a-1})$  where $gcd(v(n),e(n))=1$. 
The following result is an extension of a result used in the proof of Hilbert's tenth problem  [8,9]) to well-known general recurrences.\\
Theorem 9: $e(n)^{2}$ divides $e(nm)$ if and only if $e(n)$ divides m\\
Proof:  From the result in section 4 we can easily observe that 
since gcd(v(n),e(n))=1, $e(n)^{2}$ divides e(an) if and only if e(n) divides a. Hence the result follows\\

The following is a extension to general recurrences of th result due to Hoggatt et al. [16]\\
Theorem  10: $e(n)^{k+1}$ divides $ e(ne(n)^{k})$\\
Proof: The proof is by induction on k. The case k=1 holds because of Theorem1. Now say that   $e(n)^{k}$ divides $ e(ne(n)^{k-1})$ is true.\\
We can write $e(ne(n)^{k})=2^{1-a}e(ne(n)^{k-1})(Ke(ne(n)^{k-1})^{2}+e(n)v(n)^{e(n)-1})$ by substituting $n=ne(ne(n)^{k-1})$ 
 and $a=e(n)$ in the result from section 4. Since $e(n)^{k}$ divides $ e(ne(n)^{k-1})$ it follows that  $e(n)^{k+1}$ divides $e(ne(n)^{k-1})(Ke(ne(n)^{k-1})^{2}+e(n)v(n)^{e(n)-1})$ and hence the theorem.\\
\\
The following well known result for the Fibonacci sequence also extends naturally to a general recurrence starting in $\{0,1\}$ (also see [19] for a different approach)\\
Theorem 11: If the sequence E has starting values of ${0, 1}$, then E(a) divides E(b) if and only if a divides b. \\
Proof:- From Theorem 1, we have that for a prime p, the sequence of indices for which $E(a)=0(mod$ $p)$ form an arithmetic progression. Since E(0)=0 and a divides b then a and b must be members of the same arithmetic progession. Hence p also divides E(b). This is true for all primes that divide E(a). Hence E(a) must divide E(b). For the converse, that is a divides b, that E(a) divides E(b) can be proved by using Binet's formula. $E(a)=(p^{a}-q^{a})/k$ and likewise for E(b). So if E(a) divides E(b) then a must divide b. \\

\section{ Conclusions}

Practical applications of the Fibonacci and Lucas sequences are well-known in diverse domains including many in computer science,  number theory (including primality testing), etc.. In fact the generalized sequences are as fascinating if not more (as pointed out in [19]),  and perhaps more potential impact. What we have done in this paper is to try to explore some of the harder areas like prime power divisibility, which could have implications for many open problems in these areas.

\end{document}